\documentclass{article}

\usepackage{amsmath}
\usepackage{amsfonts}
\usepackage{amssymb}
\usepackage{theorem}
\usepackage{psfrag}
\usepackage{graphicx}

\newcommand\Rset{\mathbb{R}}
\newcommand\Nset{\mathbb{N}}

\newcommand{\Umnr}{{Up^+(n,\Rset)}}

\newcommand{\SOnr}{{SO(n,\Rset)}}

\newcommand{\OL}{{\mathcal{ O}_\Lambda}}

\newcommand{\Anr}{{Sk(n,\Rset)}}
\newcommand{\Unr}{{Up(n,\Rset)}}
\newcommand{\FS}{{\mathcal{ F}_S}}
\newcommand{\PL}{{\mathcal{ P}_{\Lambda}}}
\newcommand{\PS}{{\mathcal{ P}_{S}}}
\newcommand{\diag}{\hbox{\rm diag}\,}
\newcommand{\pia}{{\Pi_{\mathcal A}\,}}
\newcommand\GS[1]{{\big[ {#1} \big]_Q}}
\newcommand\GST[1]{{\big[ {#1} \big]_Q^T}}
\newcommand\TS[1]{{\big[ {#1} \big]_R}}
\newcommand\TSI[1]{{\big[ {#1} \big]_R^{-1}}}
\newcommand{\tr}{\hbox{tr}\,}
\newcommand{\ttI}{\mathcal{I}}
\newcommand{\ttJ}{\mathcal{J}}
\newcommand{\sgs}{Q_\ttI}
\newcommand{\sts}{R_\ttI}
\newcommand{\e}{\hbox{\rm e}\,}

{\theorembodyfont{\rmfamily}\newtheorem{prop}{Proposition}}
{\theorembodyfont{\rmfamily}\newtheorem{thm}{Theorem}}
{\theorembodyfont{\rmfamily}\newtheorem{lem}{Lemma}}
{\theorembodyfont{\rmfamily}}
{\theorembodyfont{\rmfamily}}

\newcommand{\qede}{\unskip\nobreak
\penalty50\hskip1.75em\null\nobreak\emptysquare\hfill\vrule width0pt
{\parfillskip=0pt \finalhyphendemerits=0 \par}}

\newcommand{\emptysquare}{{\hbox{\vrule height6pt width0.6pt depth0pt%
\vbox{\hrule height0.6pt width4.8pt depth0pt%
\vglue4.8pt%
\hrule height0.6pt width4.8pt depth0pt}%
\vrule height6pt width0.6pt depth0pt}}}

\begin{document}

\title{Parametrization by polytopes of intersections of orbits by conjugation}

\author{R. S. Leite \and C. Tomei}
%\footnote{Partially supported by CNPq, Faperj and Finep, Brazil} 

\maketitle

\begin{abstract}
Let $S$ be an $n\times n$ real symmetric matrix with spectral decomposition
$S\, =\, Q^T\, \Lambda\, Q$, where $Q$ is an orthogonal matrix and $\Lambda$ is diagonal with 
simple spectrum $\{\lambda_1,\ldots, \lambda_n\}$. Also let ${\mathcal O}_S$ e ${\mathcal R}_S$ be
the orbits by conjugation of $S$ by, respectively, orthogonal matrices and upper triangular
matrices with positive diagonal. Denote by ${\mathcal F}_S$ 
the intersection ${\mathcal O}_S \cap {\mathcal R}_S$. We show that the map
$F\colon  {\bar{\mathcal F}_S}  \to \Rset^n$ taking 
$S' = (Q')^T\, \Lambda\, Q'$ to  $\hbox{diag}\, (Q'\, \Lambda\, (Q')^T)$ 
is a smooth bijection onto its range ${\mathcal P}_S$, the convex hull of some subset of
the $n!$ points
$\{(\lambda_{\pi(1)}, \lambda_{\pi(2)}, \ldots, \lambda_{\pi(n)}),\, \pi 
\hbox{ is a permutation}\}.$
We also find necessary and sufficient conditions for ${\mathcal P}_S$ to have $n!$ vertices. 
\end{abstract}
\bigskip

\textit{Keywords:} Integrable systems,  Toda flows, QR decomposition
\medskip

\textit{Mathematical subject classification:} Primary 58F07, Secondary 15A23

%%%%%%%%%%%%%%%%%%%%%%%%%%%%%%%%%%%%%%%%%%%%%%%%%%%%%%%%%%%%%%%%%%%%%%%%%%%%%%%%%%%%%%%%%%%

%This file contains the tags for \psfrag

\psfrag{diag124}
{\tiny
$\begin{pmatrix}
1 & 0 & 0 \cr
0 & 2 & 0 \cr
0 & 0 & 4 \cr
\end{pmatrix}$}

\psfrag{diag142}
{\tiny
$\begin{pmatrix}
1 & 0 & 0 \cr
0 & 4 & 0 \cr
0 & 0 & 2 \cr
\end{pmatrix}$}

\psfrag{diag412}
{\tiny
$\begin{pmatrix}
4 & 0 & 0 \cr
0 & 1 & 0 \cr
0 & 0 & 2 \cr
\end{pmatrix}$}

\psfrag{diag421}
{\tiny
$\begin{pmatrix}
4 & 0 & 0 \cr
0 & 2 & 0 \cr
0 & 0 & 1 \cr
\end{pmatrix}$}

\psfrag{diag214}
{\tiny
$\begin{pmatrix}
2 & 0 & 0 \cr
0 & 1 & 0 \cr
0 & 0 & 4 \cr
\end{pmatrix}$}

\psfrag{diag241}
{\tiny
$\begin{pmatrix}
2 & 0 & 0 \cr
0 & 4 & 0 \cr
0 & 0 & 1 \cr
\end{pmatrix}$}

\psfrag{m1**}
{\tiny
$\begin{pmatrix}
1 & 0 & 0 \cr
0 & * & * \cr
0 & * & * \cr
\end{pmatrix}$}

\psfrag{m2**}
{\tiny
$\begin{pmatrix}
2 & 0 & 0 \cr
0 & * & * \cr
0 & * & * \cr
\end{pmatrix}$}

\psfrag{m4**}
{\tiny
$\begin{pmatrix}
4 & 0 & 0 \cr
0 & * & * \cr
0 & * & * \cr
\end{pmatrix}$}

\psfrag{m**1}
{\tiny
$\begin{pmatrix}
* & * & 0 \cr
* & * & 0 \cr
0 & 0 & 1 \cr
\end{pmatrix}$}

\psfrag{m**2}
{\tiny
$\begin{pmatrix}
* & * & 0 \cr
* & * & 0 \cr
0 & 0 & 2 \cr
\end{pmatrix}$}

\psfrag{m**4}
{\tiny
$\begin{pmatrix}
* & * & 0 \cr
* & * & 0 \cr
0 & 0 & 4 \cr
\end{pmatrix}$}

\psfrag{m*1*}
{\tiny
$\begin{pmatrix}
* & 0 & * \cr
0 & 1 & 0 \cr
* & 0 & * \cr
\end{pmatrix}$}

\section{Introduction}

In the seventies, interest in the Toda lattice 
--- a Hamiltonian system of equations for 
the motion of $n$ particles in the line ---
led researchers to consider in detail the geometry of Jacobi matrices 
(real, symmetric tridiagonal matrices
with strictly positive elements immediately below and above the diagonal).
The connection was made by Flaschka \cite{[F]}, which constructed a change of variables 
from physical phase space $\Rset^{2n}$ to Jacobi matrices, converting 
the Toda equations into an evolution with the remarkable property that the 
eigenvalues of the matrices along an orbit are constant. 
After handling the motion of the center of mass of the particle system,
the change of variables transfers the
standard symplectic structure on $\Rset^{2n-2}$  to submanifolds of $n \times n$
Jacobi matrices with trace zero. By proving that (symmetric
functions of the) eigenvalues taken as Hamiltonians induce commuting flows, 
Flaschka showed that the particle system is completely integrable on both
the original variables and on Jacobi matrices. On the matrix phase space,
the Liouville-Arnold theorem for integrable systems 
(which still holds, despite the noncompactness of the invariant tori)
becomes the statement that the set ${\mathcal J}_\Lambda$ of $n \times n$ Jacobi matrices 
with simple spectrum equal to the spectrum of a diagonal matrix 
$\Lambda = \diag (\lambda_1 > \lambda_2 >\ldots> \lambda_n)$ 
is diffeomorphic to $\Rset^{n-1}$. This domain of parameters for ${\mathcal J}_\Lambda$ 
may be identified with the evolution times of $n-1$ commuting Toda flows starting from
a fixed initial condition: there is thus a natural {\it Toda action} on ${\mathcal J}_\Lambda$. 
Later, Moser \cite{[M]} provided a bijective 
parametrization of ${\mathcal J}_\Lambda$ by introducing appropriate 
spectral data on Jacobi matrices. 

Further study \cite{[T]} of ${\mathcal T_\Lambda}$, the set of {\textsl{ all}}
tridiagonal matrices with fixed simple spectrum $\Lambda$,  
required a detailed understanding of the closure $\bar{{\mathcal J}_\Lambda}$
in the vector space of real, symmetric, tridiagonal matrices.
By combining  topological and combinatorial arguments, $\bar{{\mathcal J}_\Lambda}$
was shown to be homeomorphic to the polytope $\PL$, the convex span of the $n!$ 
points in $\Rset^n$ of the form
$(\lambda_{\pi(1)}, \ldots, \lambda_{\pi(n)}),$
where $\pi$ is an arbitrary permutation of the set $\{1,\ldots,n\}$.

The appearance of the convex polytope $\PL$ suggested that this result should
be interpreted in terms of symplectic torus actions, following the
fundamental papers of Guillemin and Sternberg \cite{[GS]} and Atiyah \cite{[A]}. 
This was accomplished by Bloch, Flaschka and Ratiu \cite{[BFR]}, 
which presented a beautiful explicit correspondence between $\bar{{\mathcal J}_\Lambda}$
and $\PL$. In order to describe it, we set some notation.
Take $S$ real symmetric with the same (simple) spectrum as $\Lambda$
and consider the spectral decomposition $ S = Q^T\, \Lambda\, Q,$ where $Q \in \SOnr,$
the group of orthogonal matrices with determinant equal to 1.
For another decomposition 
$ S\, =\, \hat Q^T\, \Lambda\, \hat Q,\quad \hat Q  \in \SOnr, $
we must have
$Q = \hat Q\, D$, where $D$ is a diagonal matrix with diagonal entries 
equal to either 1 or -1. Thus, the matrix $\tilde S = Q^T\, \Lambda\, Q$ 
is actually dependent on the choice of $Q$, but its diagonal entries are not. 
Also, the trace of $\tilde S$ is the trace of $\Lambda$. Let ${\mathcal B}_\Lambda$ be the 
hyperplane in $\Rset^n$ of vectors whose coordinates add to $\tr \Lambda$. On the manifold
$\OL$ of real symmetric matrices with spectrum equal to $\Lambda$,
we define the \textsl {BFR map}
$$F_{\OL} \colon {\OL} \to {\mathcal B}_\Lambda 
\quad \hbox{ given by } \quad F(S) = \diag \tilde{S}.$$
Clearly, $F_{\OL}$ is a smooth map.

\begin{thm}{\textrm{ (Bloch-Flaschka-Ratiu)}}  
The restriction $F_{\bar{{\mathcal J}_\Lambda}}\colon \bar{{\mathcal J}_\Lambda} \to \PL$ 
is a homeomorphism which is a diffeomorphism between the interiors of domain and range.
\end{thm}

The key point in the proof of \cite{[BFR]} is the fact that, 
on an appropriate K\"ahler manifold, 
a twisted version of the Toda action is the holomorphic continuation 
of the standard diagonal torus action. The result then follows from 
Atiyah's convexity theorem on extensions of moment maps on K\"ahler 
manifolds \cite{[A]}. In principle, different restrictions of $F$ 
could also be shown to be diffeomorphisms by the arguments in \cite{[BFR]}, 
but the issue was not considered.

The BFR map brings to mind the simpler Schur-Horn map, 
$G: \OL \to {\mathcal B}_\Lambda$, taking $S$ to $\diag S$ --- 
the fact that this map is surjective is the classic
Schur-Horn theorem ([S23], [H54]). When restricted to $\bar{{\mathcal J}_\Lambda}$,
this map is far from being either injective or surjective.

Closer to the spirit of the BFR theorem, the object of this paper is to provide an 
elementary proof that, instead, many restrictions of the BFR map are diffeomorphisms, 
taking values on appropriate convex polytopes: a special case yields
the BFR theorem. 
Take $S \in {\OL}$. Consider the {\textsl{ slice}} $\FS$ through $S$,
\begin{equation}\label{intersecao}
\{\, Q^T\, S\, Q,\, Q\in\SOnr \,\} \cap \{\, R\, S\, R^{-1},\, R\in\Umnr\,\},
\end{equation}
where $\Umnr$ is the group of upper triangular matrices with positive diagonal entries.
Consider  the \textsl{ accessible vertices}, which are
the diagonal matrices belonging to the closure of a slice $\bar \FS$.
The image under $F_{\OL}$ of an accessible vertex is an \textsl{ extremal vertex}. 
Finally, let
$\PS$, the \textsl{spectral polytope of } $S$, be the convex hull of the extremal
vertices.

\begin{thm} Let $S$ be a real, symmetric matrix, with simple spectrum. The restriction
$F_{\bar{\FS}}\colon{\bar{\FS}} \to \PS$ is a homeomorphism,
which is a diffeomorphism between the interiors of domain and range.
\end{thm}

A special case of this result
occurs for \textsl{spectrally complete} matrices in $\OL$: these are matrices 
$S = Q^T\, \Lambda\, Q,\, Q \in \SOnr,$
so that every minor of $Q$
obtained by intersecting arbitrary $k$ rows with the first $k$ columns
has nonzero determinant.

\begin{thm} Let $S$ be a spectrally complete matrix. The restriction of the BFR map
$F_{\bar{\FS}}\colon{\bar{\FS}} \to {\mathcal P}_\Lambda$ is a homeomorphism,
which is a diffeomorphism between the interiors of domain and range.
\end{thm}

Jacobi matrices are spectrally complete, as will be shown in Proposition \ref{jacobi}.
Thus, Theorem 3 indeed generalizes Theorem 1. 

Up to choices of sign, the definition of the slice through $S$ is strictly algebraic.
Still, slices have an obvious affinity with the Toda flows: $\FS$ is the set of matrices
reached by Toda flows starting from $S$ --- a precise statement,
combined with the proof of this fact, will be given in Section 2. It is here that
simplicity of the spectrum of $\Lambda$ plays its role.
We then proceed in Section 3 to the
study of the asymptotic behavior of the Toda action starting from $S$: 
the elements in the closure of $\FS$ are interpreted as equivalence classes of limits to infinity. 
From this follows the description of the
stratification into faces and subfaces of the closure $\bar \FS$. 
At this point our arguments could have taken a different route: the
asymptotic behavior yields a description of $\bar \FS$ which 
is independent of any explicit identification between $\bar \FS$ and some 
(``dressed'') reference set. 
This approach will be presented in a forthcoming paper, and will allow us to consider
scenarios in which there seems to be no natural source of convexity to 
guide our steps. In Section 4,
we describe the spectral polytope $\PS$ by characterizing both its vertices and its faces. 
It turns out that $\PS$ is obtained by chopping ${\mathcal P}_\Lambda$ by
hyperplanes parallel to faces of ${\mathcal P}_\Lambda$, but such hyperplanes are special,
in the sense that the chopping does not introduce new vertices. 
We then complete the proof of Theorem 2.

Readers familiar with the proof of the convexity theorem for the image of moment maps of
symplectic torus actions in \cite{[GS]} know that there are two basic facts to confront: 
\begin{itemize}
\item{The image of the critical set of the moment map is contained in a collection of 
hyperplanes,}
\item{The image of the map ``has no holes''.}
\end{itemize}
The fact that the map of interest {\sl is} a moment map for a symplectic torus action
is what permits the verification of both statements in the general case. 
Here  instead, we shall see from a computation in Section 2 that $F_{\OL}$ ceases to
be a diffeomorphism exactly at \textsl{reducible} matrices $S  \in {\OL}$,
i.e., matrices  which admit a proper
invariant subspace generated by vectors in the canonical basis, and
the restrictions of the BFR map $F_{{\OL}}$ to slices $\FS$ have no critical set.
The image of the boundary
of each slice $\FS$ is shown in Proposition \ref{hiperplanos} to belong 
to a very explicit set of hyperplanes.
The second step, which is equivalent to showing that there are no
images of critical components in the interior of the putative range, is a direct
consequence of the characterizaton of $\PS$ by faces. Closest in spirit to the
proof in \cite{[GS]} is the construction in Section 4 of special Lyapunov functions 
for Toda flows, the \textsl{partial traces}.

Arguments of this more elementary sort were already used in
\cite{[LRT]} for the study of specific moment maps. 
This paper is an extended version of \cite{[Le]}.

%%%%%%%%%%%%%%%%%%%%%%%%%%%%%%%%%%%%%%%%%%%%%%%%%%%%%%%%%%%%%%%%%%%%%%%%%%%%%%%%%%%%%%%%%%%

\section{The slice by a matrix}

Throughout the text, $S$ denotes an $n\times n$ real symmetric matrix with simple spectrum,
admitting the spectral decomposition $S = Q^T\, \Lambda\, Q$, where $Q$ lies in $\SOnr,$ and
$\Lambda = \diag(\lambda_1 > \ldots >\lambda_n)$.
We consider the \textsl{slice} through $S$, 
$\FS = \{\, Q^T\, S\, Q,\, Q\in\SOnr \,\} \cap \{\, R\, S\, R^{-1},\, R\in\Umnr\,\}.$
Clearly, if $S' \in \FS$ then ${\mathcal F}_{S'} = \FS$.

A nontrivial subspace of $\Rset^n$ is \textsl{canonical} 
if it contains a basis  given by a (proper) subset
of the set of canonical vectors $\{ e_1,\ldots, e_n \}$. A matrix with no invariant canonical
subspaces is \textsl{irreducible}. 
Recall that, if $S$ is a symmetric matrix with simple spectrum and $f$ is a continuous function, 
then $f(S) = p(S),$ 
where $p$ is any polynomial coinciding with $f$ on the spectrum of $S$. 

\begin{prop}\label{irredutivel}
Let $S$ be a real symmetric matrix with simple spectrum. Then $S$ is irreducible if and only if 
the only diagonal matrices which are functions of $S$ are multiples of the identity matrix.
If $S$ is irreducible, then all matrices in $\FS$ are also irreducible.
\end{prop}
\noindent \textbf{Proof}
Clearly, if $f(S)$ is diagonal with two distinct eigenvalues, it must have an invariant canonical
subspace, which is also invariant under $S$, showing reducibility. Conversely, if $S$ is reducible,
consider an invariant canonical subspace $V$ and choose a function which sends the eigenvalues of
the restriction of $S$ to $V$ (resp. $V^\perp$) to 0 (resp. to 1), giving rise to a diagonal 
function of $S$ with eigenvalues 0 and 1. For the second statement,
let $S' = R\, S\, R^{-1}\in \FS$ where $R$ is upper triangular. 
Suppose that $S$ is irreducible and 
that, by contradiction, $S'$ is reducible. Let $p$ be a polynomial such that $p(S')$ is a 
diagonal matrix which is not a multiple of the identity matrix $I$. Since $S = R^{-1}\, S'\, R$, 
$$
p(S)\, =\, p(R^{-1}\, S'\, R)\, =\, R^{-1}\, p(S')\, R,
$$
so that $p(S)$ is a symmetric, upper triangular matrix, and hence, diagonal. The spectra of 
$p(S)$ and $p(S')$ are the same, which shows that $p(S)$ is not a multiple of $I$ --- a 
contradiction.
\qede \bigskip

Let $S$ be an invertible symmetric matrix. Consider the (unique) QR decomposition of $S$ as
the product of two matrices, one in $\SOnr$ and one in $\Umnr$: $S = Q\, R$. We denote the factors
$Q$ and $R$ by, respectively, $\GS{S}$ and $\TS{S}$, and its inverses by $\GST{S}$ and $\TSI{S}$. 
With an appropriate choice of signs, the matrix $\GS{S}$ is the matrix obtained by applying 
the Gram-Schmidt process to the columns of $S$. The QR factors depend smoothly on
the invertible matrices $S$.

\begin{prop} \label{fatia}
Let $S$ be an $n\times n$ real symmetric irreducible matrix with simple spectrum. Then
$$\FS  =
\{ \GST{ \e^{p(S) } }\, S\, \GS{ \e^{ p(S) }  },\, \hbox{\rm for polynomials } 
p \colon \Rset \to \Rset\} .$$
Moreover, $\GST{ \e^{p(S)} }\, S\, \GS{ \e^{p(S) } }=  
\GST{ \e^{{\tilde p}(S)} }\, S\, \GS{ \e^{ {\tilde p}(S) } }$
if and only if $p = \tilde p + b$, for some $b \in \Rset$. 
\end{prop}

\noindent \textbf{Proof}
Let $S'\in\FS$: there are matrices $Q\in\SOnr$ and $R\in\Umnr$ such that
$S'= Q^T\, S\, Q = R\, S\, R^{-1}$. Then $S\, Q\, R = Q\, R\, S$ and, since $S$ has simple
spectrum, $QR = f(S)$ for some polynomial $f$ with $f(S)$ invertible. Thus,
$Q = \GS{ f(S) }$ and $R =  \TS{ f(S) }$ and 
$$
\FS \subset \{ \GST{ f(S) }\, S\, \GS{ f(S) },\,  \hbox{\rm for polynomials } f \}.
$$

Suppose that $h = \theta\, f$ where $\theta$ take the values $\pm\, 1$ on the spectrum of $S$.
Then $\GST{ f(S) }\, S\, \GS{ f(S) }\, =\, \GST{ h(S) }\, S\, \GS{ h(S) }$,
since $\GS{ h(S) }\, =\, \theta(S)\, \GS{ f(S) }$ and $\theta(S)$ is orthogonal.
The first claim of the proposition follows: $f$ and $|f|$ 
obtain the same matrix in $\FS$.
Now take $f$ and ${\tilde f}$ positive on the spectrum of $S$ so that
$$
\GST{ f(S) }\, S\, \GS{ f(S) }\, =\, \GST{ {\tilde f}(S) }\, S\, \GS{ {\tilde f}(S) }
$$
and, from the computations in the beginning of the proof,
$$
\TS{ f(S) }\, S\, \TSI{ f(S) }\, =\, \TS{ {\tilde f}(S) }\, S\, \TSI{ {\tilde f}(S) }.
$$
The first identity yields
$\GS{ f(S) }\, \GST{ {\tilde f}(S) }\, S\, =\, S\, \GS{ f(S) }\, \GST{ {\tilde f}(S) }$, 
which shows that
$\GS{ f(S) }\, \GST{ {\tilde f}(S) }$ is a polynomial $k(S)$:
$\GS{ f(S) }\, =\, k(S)\, \GS{ {\tilde f}(S) }.$
Thus $k(S)$ is orthogonal and symmetric, and its spectrum is contained in $\{-1, +1\}$. 

By the second identity,
$S\, \TSI{ f(S) }\, \TS{ {\tilde f}(S) }\, =\, \TSI{ f(S) }\, \TS{ {\tilde f}(S) }\, S$.
The matrix $E = \TSI{ f(S) }\, \TS{ {\tilde f}(S) }\in\Umnr$ is a polynomial of $S$ 
and hence is symmetric.
Thus $E$ is actually positive diagonal and 
$\TS{ f(S) }\, E\,  =\, \TS{ {\tilde f}(S) }$.
A simple algebra obtains $f(S)\, E\, =\, k(S)\, {\tilde f}(S)$,
so
$$
E\, =\, (f^{-1}\circ k\circ  {\tilde f})(S).
$$
Since $E$ is a positive diagonal polynomial of the irreducible matrix $S$, we must have
$E = a\, I$, $a > 0$. Comparing spectra of $E = a\, I = (f^{-1}\circ k\circ {\tilde f})(S)$, 
we see that $k(S)$ is the identity, and, finally, $a\, f = {\tilde f}$.
By taking logs, convert the statements about $f$'s in statements about $e^p$'s.
\qede \bigskip

Consider the additive homomorphism $\tau \mapsto p_{\tau}$, taking a vector $\tau \in \Rset^n$
to a polynomial $p_{\tau}$ with values $\tau_1,\ldots, \tau_n$ 
on the spectrum $\lambda_1,\ldots, \lambda_n$. The quotient of the vector space of
polynomials of degree $n-1$ by constants will be represented by
$\Rset^n_0 = \{(\tau_1,\ldots,\tau_n) \in \Rset^n, \sum_i \tau_i = 0 \}$.
Consider now the map 
\begin{align*}
\Phi  \colon &(\Rset^n_0,+) \times {\OL}   \to  {\OL} \\
             &(\tau, S)  \mapsto \GST{ \e^{p_\tau(S) } }\, S\,\GS{ \e^{p_\tau(S) } }. \\
\end{align*}
This map, the \textsl{Toda action}, describes the flows associated to the 
Hamiltonians of the standard Toda hierarchy. for the convenience of the reader, 
we sketch the basic facts about Toda flows: an excellent reference is \cite{[Sy80]}.
Consider the decomposition
$S = S_u +  S_d + S_\ell$, where $S_u$, $S_d$ and $S_\ell$ are, respectively, the strictly upper 
triangular, diagonal and strictly lower triangular parts of $S$. Define the projections
$$
{\Pi_{\mathcal A}\,} S\, =\, S_\ell - S_\ell^T\in\Anr, \quad
{\Pi_{\mathcal U}\,} S\, =\, S_u +S_d + S_\ell^T\in\Unr,
$$
where $\Anr$ and $\Unr$ are the real vector spaces of skew-symmetric and upper triangular
matrices respectively. 
The solution at $t$ of the ordinary differential equation  
\begin{equation}\label{toda}
\frac{\hbox{d}}{\hbox{d}\, t}\,  S(t)  = [\, S(t)\, ,\,  {\Pi_{\mathcal A}\,} p(S(t))\, ], \quad
S(0)  = S\, .
\end{equation}
is  known to be equal to 
$$
S(t)\, =\, \GST{ \e^{ t\, p(S) } }\, S \, \GS{ \e^{ t\, p(S) } }\,
=\, \TS{e ^{ t\, p(S) } }\, S\, \TSI{ \e^{ t\, p(S) } }.
$$

The next proposition collects some obvious consequences of this formula.

\begin{prop}\label{imersao}
The matrix $\Phi(\tau, S)$ equals the solution at $t=1$ of the differential equation
above where $p = p_{\tau}$.
The slice $\FS$ is the set of all matrices reached by a Toda flow starting at $S$.
Toda flows commute: for $\tau_1,\tau_2\in\Rset^n_0$,
$$
\Phi(\tau_1,\, \Phi(\tau_2,\, S))\, =\, \Phi(\tau_1\, +\, \tau_2\, S).
$$
In particular, $\Phi$ is a group action on ${\OL}$ which preserves each slice $\FS$.
Each restriction $\Phi_S: \Rset^n_0 \to \FS$ is an injective immersion.
\end{prop}

\noindent \textbf{Proof}
For commuting matrices $M_1= \e^{p_{\tau_1}(S)}$ and $M_2 = \e^{p_{\tau_2}(S)}$, 
\begin{align*}
\GS{ M_2 \GS{ M_1  } } &= \GS{ M_2 \GS{ M_1 } \TS{ M_1 } }
= \GS{ M_2 M_1 }\cr 
&=\GS{ M_1 M_2 } = \GS{ M_1 \GS{ M_2  } },\cr
\end{align*}
from which commutativity of the action follows. from the previous proposition,
$\Phi$ is goes down to the quotient $\Rset^n_0$ and, for the same reason, it is
injective in $\FS$.
To see that $\Phi_S: \Rset^n_0 \to \FS$ is an immersion, 
begin by differentiating $\GS{ \e^{ p_{\tau}(S) } }$ in the 
variable $\tau$ along the direction $\rho\in\Rset^n_0$,
\begin{align*}
D_\tau\, \left( \GS{ \e^{ p_{\tau}(S) } } \right)\, (\rho) & = 
\GS{ \e^{ p_{\tau}(S) } }\, {\Pi_{\mathcal A}\,}\left( \GST{ \e^{ p_{\tau}(S) } }\,
p_{\rho}(S)\, \e^{ p_{\tau}(S) }\,
\TSI{ \e^{ p_{\tau}(S) } } \right) \\
& = \GS{ \e^{ p_{\tau}(S) } }\, {\Pi_{\mathcal A}\,}\left( \GST{ \e^{ p_{\tau}(S) } }\,
p_{\rho}(S)\,  \GS{ \e^{ p_{\tau}(S) } } \right), \\
\end{align*}
where, for $M = \e^{ p_{\tau}(S) }$, we used the identity
\begin{equation}\label{deri}
\frac{\hbox{d}}{\hbox{d}\, \tau}\, {\GS{\, M\, }}\, =\, 
\GS{\, M\, }\, {\Pi_{\mathcal A}\,} (\, \GST{\, M\, }\, (\frac{\hbox{d}}{\hbox{d}\,\tau}\, M)\, 
\TSI{\, M\, }\, ).
\end{equation}
Finally, differentiating $\Phi_{S}$ at $\tau=0$ along $\rho$,
$$
D_\tau\, \Phi_{S}|_{\tau=0} (\rho)\, =\, [ {\Pi_{\mathcal A}\,} p_{\rho}(S)\, ,\, S ].
$$
The map $\Phi_S$ is a immersion if and only if $[{\Pi_{\mathcal A}\,} p_{\rho}(S), S ]\, =\, 0$:
the skew-symmetric matrix ${\Pi_{\mathcal A}\,}  p_{\rho}(S) $ is thus a function of $S$, 
hence symmetric,
and we must have ${\Pi_{\mathcal A}\,} p_{\rho}(S)= 0$.
So $p_{\rho}(S)$ is symmetric upper triangular, hence  diagonal. 
As $p_{\rho}(S)$ is a function of the irreducible matrix $S$, 
$p_{\rho}(S)$ must be a multiple of the identity matrix,
by Proposition \ref{irredutivel}, 
and its eigenvalues (the coordinates of the vector $s$ of $\Rset^n_0$) must be zero: $s=0$.
\qede \bigskip

%%%%%%%%%%%%%%%%%%%%%%%%%%%%%%%%%%%%%%%%%%%%%%%%%%%%%%%%%%%%%%%%%%%%%%%%%%%%%%%%%%%%%%%%%%%%%

\section{Asymptotics of the Toda action}

Given a matrix $S = Q^T\, \Lambda\, Q$, we consider 
the matrices in the closure of $\FS$. These matrices will be obtained
by taking limits of $\Phi(\tau^k,S)$ for appropriate choices of vectors $\tau^k \in \Rset^n_0$.
From the solution of the Toda flows,
\begin{align*}
\Phi(\tau^k, S) &= \GST{ e^{ p_{\tau^k} (S) } }\, S\, \GS{ e^{ p_{\tau^k} (S) } }\, 
=\,
\GST{Q^T\, e^{p_{\tau^k}(\Lambda)}\, Q}\, Q^T\, \Lambda\, Q\, 
\GS{Q^T\, e^{ p_{\tau^k}(\Lambda)}\, Q}\\
&= \GST{e^{p_{\tau^k}(\Lambda)}\,  Q}\, \Lambda\, \GS{e^{p_{\tau^k}(\Lambda)}\, Q}. \\
\end{align*}
We choose a different representative for classes in the quotient by constants
of the space of polynomials: instead of points in $\Rset^n_0$, we take vectors\\ 
$ \sigma = (\sigma_1,\ldots,\sigma_n)$ with largest coordinate equal to 0. Vectors $\tau^k$ are
then replaced by their counterparts $\sigma^k$ in the formulae above and,
by definition, $p_{\sigma^k}(\Lambda) = \diag(\sigma_1^k,\ldots,\sigma_n^k)$, and
$\GS{ e^{p_{\sigma^k}(\Lambda)}\, Q} = \GS{D^k\, Q} = Q_k$, for some positive diagonal matrix
$D^k$ with largest entry equal to 1. 

We are thus led to consider limits $Q_{\infty}$ of sequences of orthogonal matrices
$Q_{k} = \GS{ \e^{\diag(\sigma^k)}\, Q }= \GS{ D^k \, Q}$ 
where $\{D^k, k\in\Nset \}$ is a sequence of positive diagonal matrices 
$D^k = \diag (d_1^k, \ldots, d_n^n)$ with entries in $(0,1]^n$.
A sequence of $n\times n$ diagonal matrices 
is \textsl{ normalized}
if its diagonal entries belong to $(0,1]$ and each matrix has some entry equal to 1.
Matrices in  $\bar \FS - \FS$ correspond to limits of $\Phi(\tau^k,S)$ when $\tau_k \to \infty$,
which in turn is equivalent to
the asymptotic property that some diagonal entry of the sequence $\{D^k, k\in\Nset\}$ 
goes to 0. 

An \textsl{ordered partition} $\ttI = (I_0, I_1, \ldots, I_p)$ 
is a partition of $\{1,2,\ldots,n\}$ in subsets  $I_0, I_1, \ldots, I_p$, 
taken in a prescribed order. Let $\{D^k, k\in\Nset \}$ be a 
sequence of $n\times n$ diagonal matrices with diagonal entries $d^k_i$ in $(0,1]$.
Such a sequence admits an \textsl{equiasymptotic partition} if there is an ordered partition of
$\{1,2,\ldots,n\}$ such that, for $i_\alpha\in I_\alpha$ and $i_\beta\in I_\beta$, the quotient
$d^k_{i_\beta} / d^k_{i_\alpha} $ goes to zero for $k \to \infty$
%$$
%\lim_{k\to\infty} \frac{d^k_{i_\beta} }{ d^k_{i_\alpha} }
%$$
if $\alpha < \beta$ or to a nonzero real number if $\alpha = \beta$. 
Thus, two indices $\alpha$ and $\beta$ in the same subset $I_j$ 
label  diagonal positions that have \textsl{comparable
asymptotic behavior}; also, the positions indexed by $I_{j+1}$ decrease to 0 
faster than the positions indexed by $I_j$. 
Ordered partitions for which the number of subsets $p+1$ is different
from 1 will be called \textsl{proper}.

\begin{lem}\label{equi}
Let $\{ D^k, k\in\Nset\}$ be a normalized sequence of diagonal matrices.
Then $\{D^k, k\in\Nset\}$ has a convergent subsequence $\{ D^\ell, \ell\in\Nset\}$ 
admitting an equiasymptotic partition $\ttI$. The limit matrix 
may be chosen to have one diagonal entry equal to 1. Also, the limit matrix has a
diagonal entry equal to zero if and only if the partition $\ttI$ is proper.

\end{lem}
\noindent \textbf{Proof} 
If the sequence $\{ D^k, k\in\Nset\}$ admits a subsequence converging to a matrix
having diagonal entries in $(0,1]$, the result is clear and the equiasymptotic partition
has a single subset (i.e., $p=0$). Otherwise, we use induction: in particular, 
the proof will obtain proper partitions. 
Since each matrix $D^k$ has a diagonal entry equal to 1,  
there is an index $i_1$ so that $d^k_{i_1} = 1$ for infinitely many $k$'s
forming an infinite subset $K_1$ of $\Nset$. Without loss, we may take $i_1 = 1$.
Suppose by induction that we have already partitioned the first $j$ indices, by assigning
them into subsets of indices with comparable asymptotic behavior, 
and we are left with a subsequence of diagonal matrices labeled by the 
infinite set $K_j$. Also, suppose that the subsets of indices are ordered
according to decreasing asymptotic behavior.
Consider now index $j+1$. For $m \in \{1,2,\ldots,j\}$,
The quotients $\{d^k_{j+1}/d^k_m, k \in K_j\}$ must accumulate in one 
of three possibilities: 0, a nonzero real number or $\infty$. Choose a 
convergent subsequence of quotients, obtaining in the process an infinite 
subset $K_{j+1,m} \subset K_j$. If the limit is 0 (resp. $\infty$) the new 
index $j$ must belong to a subset of indices appearing after (resp. before)
the subset to which the index $m$ belongs. If the limit is
a nonzero real number, $j+1$ and $m$ belong to the same subset of indices.
After repeating this process for all indices $m = 1,2,\ldots,j$, we get to 
know which subset of indices (possibly a new one) contains the index $j+1$.
Finally, set $K_{j+1} = K_{j+1,j}$.
\qede \bigskip

A \textsl{boundary sequence} is a sequence $\{D^k, k\in\Nset\}$ of normalized 
diagonal matrices admitting a proper equiasymptotic partition. 
We will see in Theorem \ref{front} that boundary sequences indeed give rise to 
boundary points of slices.

Before analyzing in full generality the asymptotic behavior of the Toda action,
or, more geometrically, the closure of a slice, we present two examples.
In the first one, the slice is an open topological hexagon, 
as are  slices through $3\times 3$ Jacobi matrices with simple spectrum 
([T084],[BFR90]). In the second example, the slice 
is an open topological quadrilateral.

For $3 \times 3$ diagonal matrices, there are  12 kinds of equiasymptotic partitions:
there are six in which each partition subset $I_\ell, \ell=0,1,2$, has exactly one element, 
three in which $|I_0| = 1$ and $|I_1| = 2$ and three more, in which $|I_0| = 2$ and $|I_1| = 1$.
Let $A = \{a_1 < a_2 < \ldots < a_k \}$ be a subset of $\{1,2,\ldots, n\}$. 
Denote by $V_A$ the $k$-dimensional subspace of $\Rset^n$ generated by the canonical vectors
$e_{a_1}, e_{a_2}, \ldots, e_{a_k}$. 

For the first example, we take  the matrix $S = Q^T\, \Lambda\, Q$, where
\begin{equation}
\Lambda\, =\,
\begin{pmatrix}
4 & 0 & 0 \cr
0 & 2 & 0 \cr
0 & 0 & 1 \cr
\end{pmatrix}\quad
\hbox{\rm and}\quad
Q\, =\,
\begin{pmatrix}
\frac{ \sqrt{6} } {6} & -\frac{ \sqrt{2} } {2} & \frac{ \sqrt{3} } {3}  \cr
\frac{ \sqrt{6} } {3} & 0                      & -\frac{ \sqrt{3} } {3} \cr
\frac{ \sqrt{6} } {6} & \frac{ \sqrt{2} } {2}  & \frac{ \sqrt{3} } {3}   \cr
\end{pmatrix}.
\end{equation}

The boundary sequence $D^k =\diag(\frac{1}{k},1,\frac{1}{k^2})$
gives rise to the equiasymptotic partition $(I_0 = \{2\}, I_1=\{1\}, I_2=\{3\})$. 
The reader will have no difficulty in checking that
$$
Q_\infty\, =\, \lim_{k\to\infty} \GS{D^k\, Q}\,
=\,
\begin{pmatrix}
0 & 1 & 0 \cr
1 & 0 & 0 \cr
0 & 0 & 1 \cr
\end{pmatrix}.
$$
For $J_0=\{1\}, J_1=\{2\}, J_2=\{3\}$, $Q_\infty$ takes 
$V_{J_0}$, $V_{J_1}$ and $V_{J_2}$ to, respectively,
$V_{I_0}$, $V_{I_1}$ and $V_{I_2}$.
The matrix
$$
S_\infty\, =\, Q_\infty^T\, \Lambda\, Q_\infty\,
=\,
\begin{pmatrix}
2 & 0 & 0 \cr
0 & 4 & 0 \cr
0 & 0 & 1 \cr
\end{pmatrix}
$$
preserves the canonical subspaces $V_{I_0}$, $V_{I_1}$ and $V_{I_2}$. 
The slice by $S_\infty$ consists of the single matrix $Q_\infty$, which 
turns out to be a  vertex of $\bar\FS$. Also, different choices 
of boundary sequences $D^k$ yielding the same ordered partition 
$(I_0 = \{2\}, I_1=\{1\}, I_2=\{3\})$ obtain the same limit matrices 
$Q_\infty$ and $S_\infty$. Other choices of boundary sequences
yielding  singletons $I_0, I_1$ and $I_2$ obtain 
the remaining five diagonal matrices with the same spectrum as $\Lambda$, 
and complete the set of accessible vertices of $\bar \FS$. 

Consider now
the boundary sequence $D^k =\diag(\frac{a}{k}, \frac{b}{k},1)$,
for $a,b \in (0,1]$, with gives rise to the equiasymptotic partition 
$(I_0=\{3\}, I_1=\{1,2\})$. Then
$$
Q_\infty^{a,b}\, =\, 
\begin{pmatrix}
\begin{matrix} 0 \cr 0 \cr \end{matrix} \hspace*{0.1cm} 
\left[ \begin{matrix}
-a\, \frac{\sqrt{2}}{2} & a\, \frac{\sqrt{2}}{2} \cr
-b\, \frac{\sqrt{2}}{2} & -b\, \frac{\sqrt{2}}{2} \cr
\end{matrix} 
\right]_Q \cr
1\hspace*{1.1cm} 0\hspace*{1.3cm}  0\hfill  \cr
\end{pmatrix}.
$$
The matrix $Q_\infty^{a,b}$ depends only on the quotient $a/b$.
Also, for the partition $(J_0 = \{1\}, J_1 = \{2,3\})$, the matrix
$S_\infty^{a,b} = (Q_\infty^{a,b})^T\, \Lambda\, Q_\infty^{a,b}$ keeps
$V_{J_0}$ and $V_{J_1}$ invariant.
The spectra of the restrictions $S_\infty^{a,b}|_{V_{J_0}}$ 
and $\tilde S_\infty^{a,b}|_{V_{J_1}}$ are, respectively,
$\{1\}=\{\lambda_i, i\in{\tilde I}_0\}$ and $\{4,2\}=\{\lambda_i, i\in{\tilde I}_1\}$. 
The slice by $S_\infty^{a,b}$ for fixed $a$ and $b$ is a topological side of $\bar \FS$ 
joining  vertices $\diag(1,4,2)$ and $\diag(1,2,4)$. Also, this slice
is exactly the set of matrices $S_\infty^{a',b'}$, for 
$a',b' \in (0,1]$.  The other two partitions 
$(I_0=\{2\}, I_1=\{1,3\})$ and $(I_0=\{1\}, I_1=\{2,3\})$ 
correspond respectively to slices of matrices preserving the subspace generated 
by $e_1$ with spectra $\{2\}$ and $\{4\}$. 

So far, we obtained three of the six sides of the topological hexagon $\bar \FS$. 
The boundary sequence $D^k = \diag(a,b,\frac{1}{k})$ for constants $a,b \in (0,1]$,
with equiasymptotic partition $(I_0=\{1,2\}, I_1=\{3\})$ and its appropriate permutations
yield the remaining three sides. The upshot is the picture below: slices in
the boundary of $\FS$ correspond bijectively to the possible equiasymptotic partitions. 
 
\bigskip
\begin{figure}[h]
\begin{center}
\includegraphics*[width=10cm]{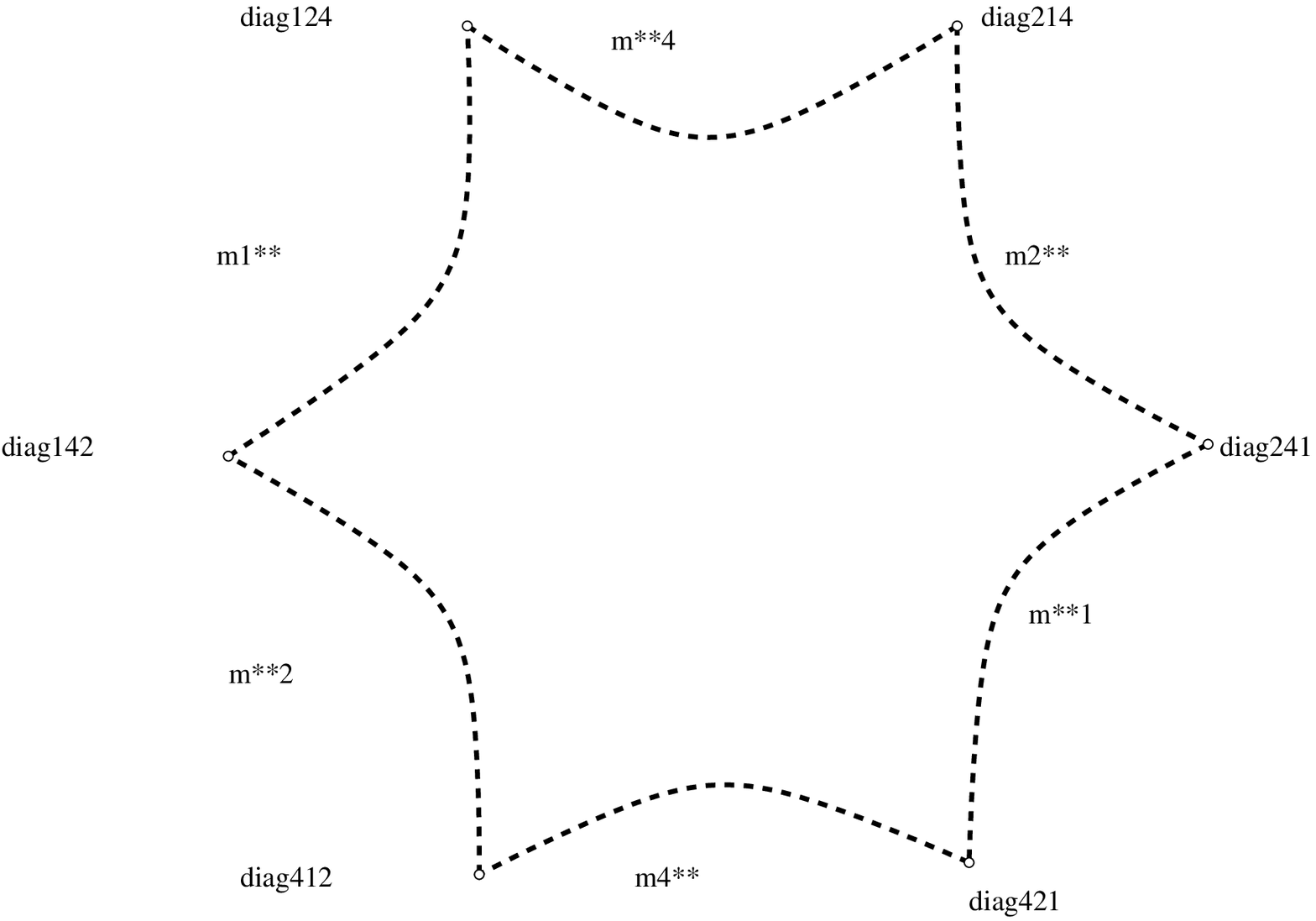}
\caption{A spectrally complete slice}
\end{center}
\end{figure}
\bigskip

For the second example, we take $S = Q^T\, \Lambda\, Q$ where
\begin{equation}
Q \, =\,
\begin{pmatrix}
\frac{\sqrt{3}}{3} & \frac{\sqrt{6}}{6}      & -\frac{\sqrt{2}}{2} \cr
\frac{\sqrt{3}}{3} & \frac{\sqrt{6}}{6}      & \frac{\sqrt{2}}{2} \cr
\frac{\sqrt{3}}{3} & -\frac{\sqrt{6}}{3}     & 0  \cr
\end{pmatrix}.
\end{equation}
Again, we begin with the boundary sequence 
$D^k = \diag(\frac{1}{k},1,\frac{1}{k^2})$,
with equiasymptotic partition $(I_0 = \{2\}, I_1 = \{1\},I_2=\{3\})$. Now
$$
Q_\infty\, =
\begin{pmatrix}
0 & 0 & 1 \cr
1 & 0 & 0 \cr
0 & 1 & 0 \cr
\end{pmatrix}, \quad
S_\infty\, =\,
\begin{pmatrix}
2 & 0 & 0 \cr
0 & 1 & 0 \cr
0 & 0 & 4 \cr
\end{pmatrix}.
$$
For the partition $(J_0=\{2\}, J_1=\{3\}, J_2=\{1\})$, the matrix $Q_\infty$ takes 
$V_{J_0}$, $V_{J_1}$ and $V_{J_2}$ to, resp., $V_{I_0}$, $V_{I_1}$ and 
$V_{J_2}$. 
Notice however that the boundary sequence $D^k=\diag(\frac{1}{k^2},1,\frac{1}{k})$, 
admitting $(I_0=\{2\}, I_1=\{3\}, I_2=\{1\})$ as equiasymptotic partition, 
also yields the same vertex $S_\infty$ of $\FS$. 
The other four choices of singleton partitions yield three more diagonal matrices: 
$(\{1\},\{2\},\{3\})$ and
$(\{1\},\{3\},\{2\})$ yield $\diag(4,1,2)$, $(\{3\},\{2\},\{1\})$ yields
$\diag(1,2,4)$ and $(\{3\},\{1\},\{2\})$ yields $\diag(1,4,2)$.

Consider now the boundary sequence $D^k=\diag(a,b,\frac{1}{k})$, with equiasymptotic
partition $(I_0=\{1,2\}, I_1=\{3\})$. Then 
$$
Q_\infty^{a,b}\, =\, \lim_{k\to \infty} \GS{D^k\, Q_2}\, =\,
\begin{pmatrix}
* & 0 & * \cr
* & 0 & * \cr
0 & 1 & 0 \cr
\end{pmatrix},
$$
where 
$$
\begin{pmatrix}
* & * \cr
* & * \cr
\end{pmatrix}\, 
=\, \begin{pmatrix}
a\, \frac{\sqrt{3}}{3} &  -a\, \frac{\sqrt{2}}{2} \cr
b\, \frac{\sqrt{3}}{3} & b\, \frac{\sqrt{2}}{2} \cr
\end{pmatrix}_Q
$$
takes $V_{J_0}$ and $V_{J_1}$ to, resp., $V_{I_0}$ and $V_{I_1}$: here $J_0=\{1,3\}$ and
$J_1=\{2\}$. The matrix $S_\infty^{a,b}=(Q_\infty^{a,b})^T\, \Lambda\, Q_\infty^{a,b}$ keeps
the subspace $V_{J_0}$ invariant with spectrum $\{4,2\}=\{\lambda_i,\, i\in I_0\}$. 
By varying $a$ and $b$ one obtains matrices in a single slice, with extrema given by the matrices 
$\diag(4,1,2)$ and $\diag(2,1,4)$ --- notice that the previous example did not have a slice joining
these two vertices. 
The boundary sequence $D^k=\diag(\frac{a}{k}, 1, \frac{b}{k})$, which one may expect
to give rise to a side, yields instead the vertex $\diag(2,1,4)$. This is also the case of 
$D^k=\diag(1, \frac{1}{k}, \frac{1}{k})$, which yields $\diag(4,1,2)$.
Adding up, $\bar \FS$ is the topological quadrilateral in Figure 2.

\bigskip
\begin{figure}[h]
\begin{center}
\includegraphics*[width=10cm]{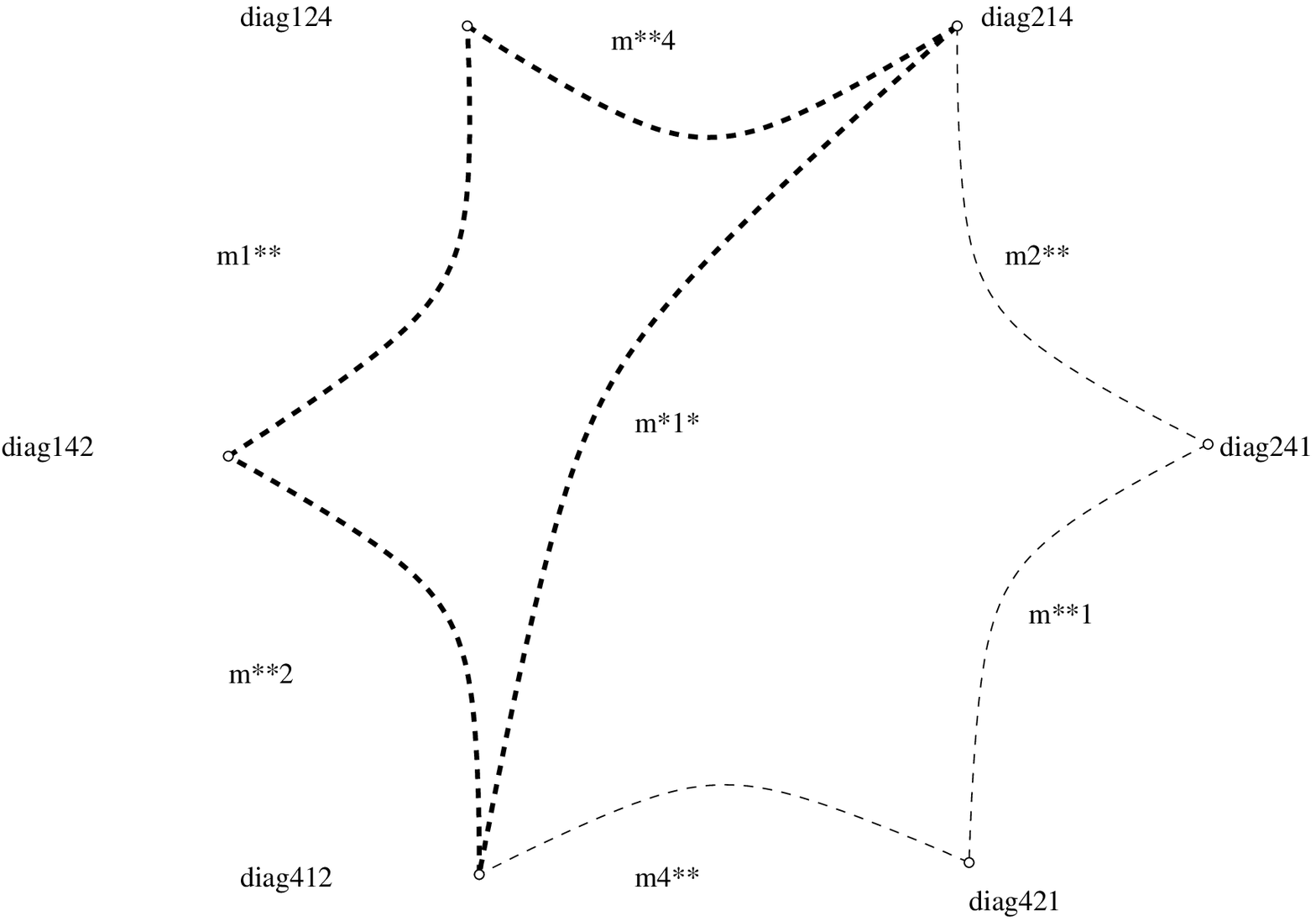}
\caption{The slice of the second example}
\end{center}
\end{figure}
\bigskip

The matrix $S$ in the first example is spectrally complete: the definition given 
in the introduction is equivalent to the more visual fact that all diagonal
matrices obtained by permuting the spectrum of $S$ belong to $\bar \FS$, as
will be proved in Proposition \ref{complete}.
The matrix in the second example is clearly not.
Theorem \ref{front} extends the pattern displayed in the examples: the
only relevant information in the computation of the limits of the Toda action 
is coded by an equiasymptotic partition and the limit values of the quotients among entries
in the same subset of the partition. 

We will make use of a slight variation of the $QR$ decomposition of an invertible $n \times n$ 
matrix $M$. For two non-empty subsets $I,J$ of $\{1, 2, \ldots, n\}$,
denote by $M_{I\times J}$ the $|I| \times |J|$  submatrix of $M$ formed by 
the elements in the rows and columns of $M$ indexed by $I$ and $J$ respectively. 
Define the \textsl{row submatrix} $M_I$ associated to $I$ as $M_I = M_{I \times \{1,\ldots,n\}}$.
Recall that $V_I$ is the span of the canonical vectors $e_i, i \in I$.
For a vector $v\in\Rset^n$, define $v_I \in V_I$ to be the orthogonal projection of $v$ in $V_I$. 
Call $m_k$ the $k^{\rm th}$ column of matrix $M$. The (unique)
$\ttI$-\textsl{sieved decomposition} 
of an invertible $n \times n$ matrix $M$ for an
ordered partition $\ttI = (I_0, I_1, \ldots, I_p)$ of $\{1,2, \ldots, n\}$ generates
\begin{enumerate}
\item{another ordered partition 
$\ttJ = (J_0, J_1, \ldots, J_p)$ with $|I_\alpha| = |J_\alpha|,$ for indices $\alpha=0,\ldots,p$,}
\item{ a matrix
$\sgs$ with orthogonal blocks $(\sgs)_{I_\alpha\times J_\alpha}$ and blocks
$(\sgs)_{I_\alpha\times (\cup_{\beta > \alpha} J_\beta)}$} equal to zero,
\item{ a matrix $\sts\in\Umnr$ such that $M\, =\, \sgs\, \sts.$}
\end{enumerate}

The factorization is obtained in steps, indexed by $\alpha= 0,\ldots, p$,
yielding invertible matrices $M_0,\ldots, M_p = \sgs$. 

Set $\alpha = 0$. Compute a matrix ${\tilde M}_0$ as follows. Start from the first column ($k = 1$).
If $(m_k)_{I_0}$, the restriction of column $k$ of $M$ to the subspace $V_{I_0}$, is nonzero,
subtract from the columns to the right of $m_k$
a multiple of $m_k$ so as to make their restriction to $I_0$ orthogonal to $(m_k)_{I_0}$. 
Repeat the procedure after increasing $k$ up to $n-1$. The resulting matrix is $\tilde{M_0}$.
Call $J_0$ the set of indices of the nonzero (orthogonal, but possibly not normal) 
columns of the row matrix $({\tilde M}_0)_{I_0}$.
Now obtain $M_0$ by dividing the columns  indexed by $J_0$ of ${\tilde M}_0$ 
by positive numbers, chosen as to make the submatrix $(M_0)_{I_0 \times J_0}$ orthogonal.
 
The fact that this computation can be performed follows from simple linear algebra:
since $M$ is invertible, the row matrix $M_{I_0}$ has row rank $|I_0|$ and hence
column rank $|I_0|$. Also, under a lexicographic ordering of $n$-uples, the set $J_0$
indexes the first subset of independent columns of $M_{I_0}$. Notice that the 
submatrix $(M_0)_{I_0 \times J_0^C}$ is zero. It is easy to see that $M_0 = M\, R_0,$
for an appropriate $R_0 \in \Umnr$.

In the second step ($\alpha=1$), a similar procedure obtains an invertible matrix $M_1$.
Let $k$ range in $J_0^C$ in increasing fashion. If for a given
$k$, the restriction $((m_0)_k)_{I_1}$ is nonzero, multiples of $(m_0)_k$ should be subtracted
from each column $k+p \in J_0^C, p>0$ to obtain vectors whose restrictions to $V_{I_1}$ are 
orthogonal to $((m_0)_k)_{I_1}$. Call the resulting matrix ${\tilde M}_1$.
Again, call $J_1$ the subset consisting of the indices of the nonzero (orthogonal) columns
of the block $({\tilde M}_1)_{I_1 \times J_0^C}$.
Now, divide by positive numbers the columns  of ${\tilde M}_1$  indexed by $J_1$
to obtain a matrix $M_1$ whose submatrix $(M_1)_{I_1\times J_1}$ is orthogonal.

The submatrix $(M_0)_{I_1\times J_0^C}$ has rank $|I_1|$ since the matrix 
$M_0$ is invertible, $\det (M_0)_{I_0\times J_0} =  \pm 1$ and $(M_0)_{I_0\times J_0^C}=0$:
the procedure above is well defined. Notice that the row matrices $(M_0)_{I_0}$ and
$(M_1)_{I_0}$ are the same. More explicitly, both $(M_0)_{I_0\times J_0^C}$ and
$(M_1)_{I_0\times J_0^C}$ are blocks of zeros, and the columns of $M_0$ indexed by $J_0$ are
left unaltered. Also, the submatrix $(M_1)_{I_1\times (J_0 \cup J_1)^C}$ is zero and
$M_1 = M_0\, R_1, R_1 \in \Umnr$. 

Iterating this process for $\alpha=2,\ldots,p$ (where on each step the range of columns to
be acted upon is, respectively, $(J_0 \cup J_1)^C, \ldots, (J_0 \cup \ldots \cup J_{p-1})^C$) 
we obtain the ordered partition 
$\ttJ = (J_0, J_1, \ldots, J_p)$ of $\{1,2,\ldots, n\}$ and the matrices $\sgs$ and 
$\sts = R_p^{-1}\, R_{p-1}^{-1}\, \ldots R_0^{-1}\in\Umnr$ with the desired properties.

\noindent \textbf{Remarks}
\begin{enumerate}
\item{By construction, any column vector $m_j, j \in J_\alpha$ has all its entries labeled by
indices in $I_0 \cup \ldots \cup I_{\alpha-1}$ equal to zero.}
\item{For an open, dense set of matrices, 
the partition $\ttJ =(J_0,J_1, \ldots,J_p)$ corresponds to the simple partition  
$n_\ttI = (\{1,2,\ldots,|I_0|\},\{|I_0|+1,\ldots,|I_0|+|I_1|\},\ldots,
\{n - |I_p|+1,\ldots,n\})$. In this case, the resulting matrix $\sgs$ is lower block triangular, 
with diagonal blocks given by orthogonal matrices.}
\item{Order the subsets of $\{1,2,\ldots, n\}$ with $k$ elements lexicographically, i.e.,
$\{i_1 < i_2 < \ldots <i_k\} < \{j_1 <j_2<\ldots <j_k\}$ if and only if 
$i_r \leq j_r$ for $r=1,\ldots,k$. The set $J_0$ is the smallest
subset of size $|I_0|$ for which the matrix $M_{I_0 \times J_0}$ is invertible. 
Similarly, $|J_1|$ is the smallest subset of size $|I_1|$ for which
$M_{(I_0 \cup I_1) \times (J_0 \cup J_1)}$ is invertible, and so on.} 
\end{enumerate}

We are ready to compute limits of the Toda action.

\begin{thm}\label{front}
Let $S = Q^T \Lambda Q$ be an $n\times n$ real symmetric irreducible matrix with simple spectrum
$\{\lambda_1>\ldots>\lambda_n\}$.
Consider a matrix $S_\infty = Q_\infty^T\, \Lambda\, Q_\infty$ in the
closure $\bar \FS$.
\begin{enumerate}
\item $Q_\infty =\lim_k \GS{D^k\, Q }$, for some normalized sequence
$\{D^k,\, k\in\Nset\}$ of diagonal matrices admitting an equiasymptotic partition
$\ttI = (I_0, I_1, \ldots, I_p)$.

\item Let $\sgs$ be the orthogonal matrix obtained by the $\ttI$-sieved  decomposition of $Q$. 
Then the matrix $Q_\infty$ splits into blocks $(Q_\infty)_{I_\alpha\times J_\beta}$,
$\alpha,\beta=0,\ldots,p$, which are equal to zero if $\alpha  \neq \beta$ and which are
orthogonal matrices for $\alpha = \beta$, given by 
$$
(Q_\infty)_{I_\alpha\times J_\alpha}\, =\, 
\lim_{k\to \infty} \GS{ D^k_{I_\alpha}\, (\sgs)_{I_\alpha\times J_\alpha} }.
$$

\item Let $\ttJ = (J_0, J_1, \ldots, J_p)$ be the
ordered partition obtained from the $\ttI$-sieved  decomposition of $Q$.
The matrix $S_{\infty}$ preserves each subspace $V_{J_\alpha}$ 
and the restriction $S_\infty|_{V_{J_\alpha}}$ 
has spectrum $\{\lambda_i,\, i\in I_\alpha\}$, $\alpha=0,1,\ldots,p$. 
\end{enumerate}
\end{thm}
\noindent \textbf{Proof} To prove item (1), combine the arguments at the beginning
of the section and Lemma \ref{equi}. Now let $Q = \sgs\, \sts$ be the $\ttI$-sieved  
decomposition of $Q$. Then
$$
Q_\infty\, =\, \lim_{k\to \infty} \GS{ D^k\, Q }\, =\, 
\lim_{k\to \infty} \GS{ D^k\, \sgs\, \sts }\,
=\, \lim_{k\to \infty} \GS{ D^k\, \sgs }.
$$

Let $i_\alpha$ be a fixed index in the set $I_\alpha$, for each $\alpha = 0, \ldots,p$.
Notice that 
$$
\lim_{k \to \infty} \frac{d^k_i}{d^k_{i_\alpha}}\, =\, \infty,\, c_i,\, 0,
$$
depending if $i \in I_\beta$, for $\beta$ less, equal or larger than $\alpha$, 
respectively.
Define $E_k$ to be a diagonal matrix so that the columns of $\sgs\, E_k$ labeled by $J_\alpha$
are the columns of $\sgs$ multiplied  by $1/i_\alpha$, for each $\alpha$.
By the presence of the blocks of zeros in $\sgs$, the matrix $Z = \lim_k D^k\, \sgs\, E^k$ 
is well defined, in the sense that all its entries converge to real numbers. 
More is true: from the convergence properties of the entries of $D^kE^k,$ 
the matrix $Z$ converges to a block matrix (with respect to the partition induced 
on rows and columns by $I_\alpha$'s and $J_\alpha$'s), whose only nonzero blocks are
the invertible submatrices $Z_{I_\alpha \times J_\alpha}$, for $\alpha=0,\ldots,p$. 
By the continuity of the Gram-Schmidt process, 
$\GS{Z} = \lim_k \GS{D^k\, \sgs\, E^k} = \lim_k \GS{D^k\, \sgs}$.
This proves (2), which in turn implies (3).
\qede \bigskip

%%%%%%%%%%%%%%%%%%%%%%%%%%%%%%%%%%%%%%%%%%%%%%%%%%%%%%%%%%%%%%%%%%%%%%%%%%%%%%%%%%%%%%%%%%%%%%%%

\section{Faces and vertices of $\bar \FS$ and its image}

From the previous theorem, a matrix $S_\infty \in \bar \FS$ is obtained from the irreducible
matrix $S = Q^T\, \Lambda\, Q$ by a limit $S_\infty = Q_\infty^T\, \Lambda\, Q_\infty$, 
where $Q_\infty = \lim_k \GS{ D^k\, Q}$, 
for some normalized diagonal sequence $\{D^k, k\in\Nset\}$.  

\begin{prop} 
With the notation above, $S_\infty$ belongs to $\bar \FS - \FS$ if and only if
the normalized sequence $\{D^k, k\in\Nset\}$ is actually a boundary sequence.
\end{prop}

\noindent \textbf{Proof}
Let $\ttI$ be the ordered partition associated to the normalized sequence $\{D^k, k\in\Nset\}$.
Suppose that $\ttI$ is proper. 
Then $S_\infty$ has at least two invariant canonical subspaces, by item (3) 
of Theorem \ref{front}, 
and thus, by Proposition \ref{fatia}, $S_\infty$ is not in $\FS$. 
Conversely, if $\ttI$ has a single subset (i.e., $p=0$), Lemma \ref{equi}
guarantees that $S_\infty$ belongs to $\FS$. 
\qede \bigskip

Thus, if $S$ is irreducible,  the matrices in $\bar \FS - \FS$ are exactly the reducible matrices.
The next result guarantees that the difference $\bar \FS - \FS$ is indeed the 
boundary $\partial \FS$, where $\FS$ 
is given the induced topology from ${\OL}$.

\begin{prop}\label{variedade}
Let $S = Q^T\,\Lambda\, Q$ be a symmetric irreducible $n \times n$ matrix
with simple spectrum. The slice 
$\FS$ is a connected $(n-1)$-dimensional manifold properly immersed in ${\OL}$. 
\end{prop}

\noindent \textbf{Proof}
To see that the slice $\FS$ is a  manifold of dimension $n-1$ with the topology
induced by ${\OL}$, it suffices to show
that $\iota \circ \Phi_S:\Rset^n_0 \to \FS \hookrightarrow {\OL}$ 
is a proper injective immersion. Because of
Proposition \ref{imersao}, we are left with showing properness. 
From Proposition \ref{fatia}, if $S$ is irreducible, all matrices in
$\FS$ are also irreducible. Let $K$ be a compact set in ${\OL}$, 
and ${\mathcal C}\subset \Rset^n_0 $ be its preimage by $\iota \circ\Phi_S$:
in particular, all matrices in $K$ are irreducible. 
As $\iota\circ\Phi_S$ is continuous, ${\mathcal C}$ is a closed set.
Now suppose that there exists 
a sequence $\{ \tau^k, k\in\Nset\}$ in ${\mathcal C}$ 
with $\lim_{k\to \infty} || \tau^k || = \infty$. The corresponding normalized sequence 
of diagonal matrices $\{D^k, k\in\Nset\}$ will admit, by Lemma \ref{equi}, a convergent subsequence
whose limit has some diagonal entries equal to zero. Thus, the associated equiasymptotic
partition $\ttI$ is proper and, from the previous proposition, 
$$
S_\infty\, =\, \lim_{k\to\infty} S_k\, =\, \lim_{k\to\infty} \, \Phi(\tau^k, S)\, 
=\, \lim_{k\to\infty}\, \GST{ \e^{p_{\tau^k}(S) } }\, S\, \GS{ \e^{p_{\tau^k}(S) } }.
$$
is a reducible matrix, belonging to the closed set $K$: a contradiction. 
Thus $\iota \circ \Phi_S$ is proper.
\qede \bigskip

Thus, matrices in the boundary of $\FS$ are reducible. In particular, they belong
to a finite collection of sets $(\OL)_J$ of matrices admitting invariant canonical
subspaces $V_J$ and $V_{J^C}$, for proper subsets $J \subset \{1,2,\ldots,n\}$. 
Each $(\OL)_J$ in turn splits into  components $(\OL)_{I,J}$, 
labeled by subsets $I \subset \{1,2,\ldots,n\}$
of eigenvalues $\{ \lambda_i, i \in I\}$ of the restriction to $V_J$ of each matrix in 
$(\OL)_J$: clearly, we must have $|I| = |J|$.
Let ${\mathcal B}_{I,J}$ be the hyperplane of points $x \in \Rset^n$ satisfying 
$\sum_{i \in I} x_i = \sum_{j \in J} \lambda_j$.
For later use, define also the half-space ${\mathcal H}_{I,J}$ to be the points $x \in \Rset^n$
for which $\sum_{i \in I} x_i \leq \sum_{j \in J} \lambda_j$.
Also, recall that ${\mathcal B}_\Lambda$ is the hyperplane 
$\sum_{i \in \{1,\ldots,n\}} x_i = \sum_{j \in \{1,\ldots,n\}} \lambda_j.$
Notice that, on ${\mathcal B}_\Lambda$, the sets ${\mathcal B}_{I,J}$ 
and ${\mathcal B}_{I^C,J^C}$ coincide.

\begin{prop}\label{quebra} Let $S^r$ be a reducible matrix, with invariant
canonical subspaces $V$ and $V^C$. Denote by $S^r_V$ and $S^r_{V^C}$ the restrictions
of $S^r$ to both subspaces. Then there is a natural identification
${\mathcal F}_{S^r} = {\mathcal F}_{S^r_V} \times {\mathcal F}_{S^r_{V^C}}.$
\end{prop}

\noindent \textbf{Proof} 
Indeed, any function $f(S^r)$ admits $V$ and $V^C$ as invariant
subspaces, and thus $\GS{f(S^r)}$ also does. 
\qede \bigskip

We will make use of the proposition above when proving statements by induction on the
dimension of the matrix.

\begin{prop}\label{hiperplanos} 
The BFR map $F_{\OL}\colon \OL \to {\mathcal B}_\Lambda$ takes the sets
$(\OL)_{I,J}$ to the (finitely many) hyperplanes ${\mathcal B}_{I,J}$.
\end{prop}

\noindent \textbf{Proof} 
This follows directly from the definition of $F_{\OL}$, 
combined with the previous proposition.
\qede \bigskip

Let $S = Q^T\, \Lambda\, Q$. 
For each proper ordered partition  $\ttI = (I_0 = I, I_1 = I^C)$,
let  $\ttJ = (J_0 = J = J(I), J_1 = J^C)$ be the corresponding  partition
obtained from the $\ttI$-sieved decomposition of $Q$. Thus $Q$ induces a map $I \mapsto J(I)$ 
between subsets of the same cardinality: the corresponding hyperplanes ${\mathcal B}_{I,J(I)}$ 
and half-spaces ${\mathcal H}_{I,J(I)}$ will be denoted by ${\mathcal B}_I$ and ${\mathcal H}_I$.
Also, denote by $F_{\bar \FS}$ the restriction of $F_{\OL}$ to $\bar \FS$. 
In the next two propositions, we describe the range of $F_{\bar \FS}$ in terms of its faces.

\begin{prop}\label{faces}
Let $S = Q^T\, \Lambda\, Q$ be an irreducible, symmetric matrix 
and consider the map $I \mapsto J = J(I)$ defined above 
induced by the $\ttI$-sieved decomposition of $Q$. 
Then the boundary of $\FS$ and the closure $\bar \FS$ are taken by $F_{\bar \FS}$ 
to the collection of hyperplanes
${\mathcal B}_\Lambda \cap (\cup_{I \subset \{1,\ldots,n\}} {\mathcal B}_I)$
and to the intersections of half-spaces
${\mathcal B}_\Lambda \cap (\cap_{I \subset \{1,\ldots,n\}} {\mathcal H}_I)$
respectively.
Finally, $F_{\bar \FS}$ takes $\FS$ to $U$, 
the interior of $ {\mathcal B}_\Lambda \cap (\cap_{I \subset \{1,\ldots,n\}} {\mathcal H}_I)$, 
injectively.
\end{prop}

\noindent \textbf{Proof}
For $I$ a proper subset of $\{1,\ldots,n\}$, consider the partial trace on $x \in {\mathcal B}_\Lambda$, 
$$\hbox{tr}_I\,(x) = \sum_{i \in I} x_i.$$
We first show that the maximum of the partial trace on the image of $F_{\bar \FS}$
equals $\sum_{j\in J(I)} \lambda_j$, and it is attained by matrices in the image of 
$\partial\FS$.
Suppose that $S' = (Q')^T\, \Lambda\, Q' \in \FS$: 
for $t \in \Rset$, consider the path in $\FS$ given by $\Phi(t\, e_I, S')$, 
where $e_I = \sum_{i\in I} e_i$.
Let $E_I = \sum_{i\in I} e_i\, e_i^T$. Clearly, $p_{e_I}(S') = (Q')^T\, E_I\, Q'$ and 
$\hbox{tr}_I\, (x) = \tr E_I\, \diag(x_1,\ldots,x_n)$.
Then
\begin{align*}
\hbox{tr}_I\,(F_{\bar \FS}(\Phi(t\, e_I, S'))) 
&= \hbox{tr}_I\,(F_{\bar \FS}(\GST{e^{t\, p_{e_I}(S')}}\, S'\, \GS{e^{t\, p_{e_I}(S')}})) \cr
&= \tr E_I\, Q'\, \GS{e^{t\, p_{e_I}(S')}}\, \Lambda\, \GST{ e^{t\,p_{e_I}(S')}}\, (Q')^T\cr
&= \tr \Lambda\, \GST{ e^{t\,p_{e_I}(S')}}\, (Q')^T\,  E_I\, Q'\, \GS{e^{t\, p_{e_I}(S')}}\cr 
&= \tr \Lambda\, \GST{ e^{t\, p_{e_I}(S')}}\, p_{e_I}(S')\, \GS{e^{t\,p_{e_I}(S')}}.\cr 
\end{align*}
Using the equality (3) of Section 2, 
we take the derivative of $\hbox{tr}_I\,(F_{\bar \FS}(\Phi(t\, e_I, S'))$ for $t = 0$:
\begin{align*}
\frac{\hbox{d}}{\hbox{d}\, t}\Bigm|_{t=0}\, \hbox{tr}_I\,(F_{\bar \FS}(\Phi(t\,e_I,S')) 
&= \tr \Lambda\, [(Q')^T\, E_I\, Q'\,, \pia p_{e_I}(S')]\cr
&= \tr \Lambda\, [p_{e_I}(S'), \pia p_{e_I}(S')]\cr
&= \sum_{i=1}^n \lambda_i ( [p_{e_I}(S'), \pia p_{e_I}(S')])_{ii}\cr
&= 2 \sum_{i=1}^n \sum_{j > i} (\lambda_i - \lambda_j) (p_{e_I}(S')_{ij})^2 \geq 0. \cr
\end{align*}
Since $\lambda_1 > \lambda_2 > \ldots > \lambda_n$, this derivative is zero if and only if
$p_{e_I}(S')^2$ is a diagonal matrix. In this case, $S'$ admits proper invariant canonical 
subspaces, by Proposition \ref{irredutivel}, and thus $S'$ is not in (the interior of) $\FS$. 
Thus, along paths  $F_{\bar \FS}(\Phi(t\, e_I, S'))$ the partial trace $\hbox{tr}_I\,$ is 
strictly increasing. 
In particular, its maximal value on the range of $F_{\bar \FS}$  equals its supremum on limit values
of paths of the form $F_{\bar \FS}(\Phi(t\, e_I, S'))$, for $S'\in \FS$.
To compute it, we take the limit 
$$
S'_{\infty} = \lim_{t\to \infty} \Phi(t\, e_I, S')
= \lim_{t\to \infty} \GST{e^{t\, p_{e_I}(\Lambda)}\, Q}\, \Lambda\, \GS{e^{t\, p_{e_I}(\Lambda)}\, Q}.
$$
Now, $t\, p_{e_I}(\Lambda)$ is a path $D(t)$ of diagonal matrices, after normalization
so that each matrix has largest diagonal entry is equal to 1. This path 
is easily seen to admit the ordered partition 
$\ttI = (I_0 = I , I_1 =  I^C)$. The $\ttI$-sieved decomposition of $Q'$ yields
the associated ordered partition $\ttJ = (J_0 = J(I), J_1 = J(I)^C)$ and the matrix $Q'_{\ttI}$. 
Set $Q'_\infty = \lim_t \GS{ D(t)\, Q_{\ttI}'}$. By Theorem \ref{front}, $Q'_\infty$
takes $V_{J_i}$ to $V_{I_i}$, $i = 0,1$. Also,
$S'_{\infty} = (Q'_{\infty})^T\, \Lambda\, Q'_{\infty}$
admits the invariants subspaces $V_{J_i}, i=0,1$, 
and its restrictions to 
have spectra $\{ \lambda_k, k\in I_i\}$, for $i=0,1$.
For matrices of this form, 
$$
\hbox{tr}_I\,(F_{\bar \FS}(S'_{\infty}))\, 
=\, \tr F_{\bar \FS}(S'_{\infty})\, E_I \, =\, \sum_{j \in J(I)} \lambda_j.$$
Thus, the supremum of the partial trace is indeed achieved, and all matrices whose image
are maximal points of the partial trace $\hbox{tr}_I\,$ admit the invariant subspaces $V_J$ and $V_{J^C}$.
Furthermore,  $F_{\bar \FS}(\FS) \subset U$ and $F_{\bar \FS}(\bar \FS) \subset \bar{U}$.

We now show injectivity of $F_{\bar \FS}$ restricted to $\FS$. 
Let $S_1,S_2\in \FS$ with $F_{\bar \FS}(S_1)=F_{\bar \FS}(S_2)$. 
We may suppose, without loss, that $S=S_1$. 
By Proposition \ref{imersao}, there is a $\rho\in\Rset^n_0$ so that 
$\Phi(\rho, S_1)=S_2$. Consider the path $F_{\bar \FS}(\Phi(t\, \rho, S_1))$ along which
we take the $t$-derivative of a weighted partial trace: 
$$
\frac{\hbox{d}}{\hbox{d}\, t}|_{t=0}\,\langle F_{\bar \FS} (\Phi(t\, \rho, S_1)\, ,\, \rho \rangle\,
=\, 2 \sum_{i=1}^n \sum_{j > i} (\lambda_i - \lambda_j) (p_{\rho}(S')_{ij})^2 \geq 0,
$$
which is strictly positive in $\FS$, unless $\rho = 0$, and then $S_1 = S_2$.
\qede \bigskip

We now consider surjectivity.

\begin{prop} Let $S$ be symmetric, irreducible.
The restriction $F_{\FS}\colon \FS \to U$ is a diffeomorphism, which extends to a homeomorphism
$F_{\bar \FS}\colon \bar \FS \to \bar{U}.$ 
\end{prop}

\noindent \textbf{Proof}
We begin by proving that $F_{\FS}$ is a local diffeomorphism. 
From Proposition \ref{variedade}, the Toda action $\Phi_S$ is a diffeomorphism from $\Rset^n_0$
to $\FS$. To prove that $F_{\FS}$ is a local diffeomorphism, we only have to 
prove that the Jacobian of $F_{\FS}\circ\Phi_S$ at $\tau=0$ is invertible. 
From the definition of the Toda action,
$$
F_{\FS}\circ\Phi_S(\tau)\, =\, \diag(Q\, \GS{ \e^{ p_\tau(S) } }\, \Lambda\, \GST{ \e^{ p_\tau(S) } }\, Q^T).
$$
Differentiating at $\tau=0$ along $\rho\in\Rset^n_0 $, we obtain
$$
D_\tau (F_{\FS}\circ\Phi_{S})\Bigm|_{\tau=0}(\rho)\, =\, 
\diag(Q\, [{\Pi_{\mathcal A}\,} p_{\rho}(S), \Lambda]\, Q^T)\, .
$$

If $D_\tau (F_{\FS}\circ\phi_{S})|_{\tau=0}$ is not injective, 
there is 
$w\in \Rset^n_0$, $w\neq 0$,  so that for all $\rho\in\Rset^n_0$
$$
\langle D_\tau (F_{\FS}\circ\Phi_{S})|_{\tau=0}(\rho)\, ,\, w \rangle\, =\, 0.
$$
Let $W = \diag(w)$. Then
\begin{align*}
0\, =\, \langle D_\tau (F_{\FS}\circ\Phi_{S})|_{\tau=0}(\rho)\, ,\, w \rangle & = 
\tr D_\tau (F_{\FS}\circ\Phi_{S})|_{\tau=0}(\rho)\, W \\ 
& =  \tr Q\, [\Pi_{\mathcal A}\, p_{\rho}(S), \Lambda]\, Q^T\, W \\
& =  \tr [\Pi_{\mathcal A}\, p_{\rho}(S), \Lambda]\, p_w(S),\\
\end{align*}
for a polynomial $p_w$. Hence
$$
\tr [\Pi_{\mathcal A}\, p_{\rho}(S), \Lambda]\, p_w(S)\, =\, 0,\quad 
\hbox{\rm for all $\rho\in\Rset^n_0$. }
$$
Taking $\rho = w$ (and hence $p_\rho = p_w$), we have
$$
0\, =\, \tr [\Pi_{\mathcal A}\, p_w(S), \Lambda]\, p_w(S)\, =\,
2\, \sum_{i = i}^n \sum_{j > i} p_w(S)_{ij}^2\, (\lambda_i - \lambda_j). 
$$
By Proposition \ref{irredutivel}, since $S$ is irreducible, $p_w(S)$ is a multiple of the 
identity matrix. Also $\tr p_w(S) = \tr Q^T\, W\, Q = 0$ since $\tr W = \tr \diag(w)$ and 
$w\in\Rset^n_0$. 
Thus, $p_w(S)=0$ and $w = 0$: the Jacobian of $F_{\FS}\circ \Phi_S$ is injective at $\tau = 0$.

Since $F_{\FS}$ is an open map, points in the boundary of $F_{\FS}(\FS)$ are necessarily 
images of boundary points in $\bar \FS$. From the previous proposition, 
points in $\partial \FS$ have to go to points in $\partial U$, which are outside of $U$
by convexity. A connectivity argument then implies that $F_{\FS}$ is surjective. 
By compactness, surjectivity of $F_{\bar \FS}$ is immediate.

Next, we show that the extension $F_{\bar \FS}\colon \bar \FS \to \bar{U}$ is injective.
Suppose $S_1, S_2 \in \bar \FS$ so that $ X = F_{\bar \FS}(S_1) = F_{\bar \FS}(S_2)$. Since
$F_{\bar \FS}$ takes interior to interior injectively and boundary to boundary, 
this may only happen if $S_1, S_2 \in \partial \FS$ and $X \in \partial U$. In particular,
from the argument in the proof of the previous proposition, $S_1$ and $S_2$ are maximal
values for some partial trace $\hbox{tr}_I\,$, and hence admit common invariant canonical subspaces $V_J$
and $V_{J^C}$.
Injectivity at the boundary now follows from induction on the dimension of the matrices:
we must have $S_1 = S_2$ when restricted to $V_J$ and $V_{J^C}$.
\qede \bigskip

The diagonal matrices in $\bar \FS$  are the \textsl{accessible vertices from } $S$. 
Each accessible vertex $\Lambda_\pi$ corresponds to a permutation $\pi$ 
of the diagonal entries of $\Lambda$. 
More precisely, there is a permutation matrix $\Pi$ with entries $\Pi_{i,j} = \delta_{i, \pi(j)}$ 
so that $\Lambda_{\pi} = \Pi^T\, \Lambda\, \Pi.$

\begin{prop}\label{vertices} 
Let $S = Q^T\, \Lambda\, Q$.
A diagonal matrix $\Lambda_\pi = \Pi^T\, \Lambda\, \Pi$ is an accessible vertex of $\FS$ 
if and only if 
the minors $Q_{\{\pi(1)\},\{1\}}, Q_{\{\pi(1),\pi(2)\},\{1,2\}},\ldots, Q$ have nonzero
determinant.
\end{prop}

\noindent \textbf{Proof}
Define $D^k$ to be a diagonal matrix whose diagonal entry $\pi(i)$ equals $k^{1-i}$.
Then, if the minors $Q_{\{\pi(1)\},\{1\}}, Q_{\{\pi(1),\pi(2)\},\{1,2\}},\ldots, Q$ 
are invertible, the sequence $\{D^k, k\in\Nset\}$ has equiasyptotic partition 
$\ttI = (\{\pi(1)\},\{\pi(2)\},\ldots,\{\pi(n)\})$ and the corresponding ordered
partition $\ttJ$ 
is $(\{1\},\{2\},\ldots,\{n\})$, from the $\ttI$-sieved decomposition of $Q$. Thus
$\lim_k \GS{ D^k\, Q} = \Pi$ and $S_\infty = \Pi^T\, \Lambda\, \Pi$ belongs to $\bar \FS$.

Now, say $\Lambda_\pi = \Pi^T\, \Lambda\, \Pi$ is an accessible vertex of $\FS$. Then, 
from Theorem \ref{front}, there exists a boundary sequence $\{D^k, k\in\Nset\}$ so that
$\Lambda_\pi = Q_\infty^T\, \Lambda\, Q_\infty$ for $Q_\infty = \lim_k \GS{ D^k\, Q}$.
Clearly, $Q_\infty = E\, \Pi,$ where $E$ is a diagonal matrix with diagonal entries equal
to 1 or -1. Equating  first columns we have that $\pm\, e_{\pi(1)}$ comes from
normalizing $\lim_k D^k\, Q\, e_1$, which implies that the entry $q_{\pi(1),1}$ of $Q$
(i.e., the minor $Q_{\{\pi(1)\},\{1\}}$) has to be nonzero. Now, by keeping
track of the determinants of the minors in rows $(\pi(1),\pi(2))$ and columns $(1,2)$
of both $Q_\infty$ and $E \, \Pi$, we learn that $Q_{\{\pi(1),\pi(2)\},\{1,2\}}$ has
nonzero determinant. The proof procedes in a similar fashion.
\qede \bigskip

The \textsl{ extremal vertices} associated to $S$ are the vectors in $\Rset^n$ with 
coordinates 
$F_{\bar \FS}( \Lambda_\pi) = (\lambda_{\pi^{-1}(1)},\ldots,\lambda_{\pi^{-1}(n)}),$ 
for all accessible vertices $\Lambda_\pi \in \bar \FS$. 
The \textsl{spectral polytope} $\PS$ associated to $S$
is the convex hull of the extremal vertices of $S$. 

\begin{prop}\label{complete}
Spectrally complete matrices are 
those for which the spectral polytope is the convex set $\PL$. 
\end{prop}

\noindent \textbf{Proof}
By definition, $S = Q^T\, \Lambda\, Q$ is spectrally complete if and only if every minor of $Q$
obtained by intersecting any $k$ rows with the first $k$ columns has nonzero determinant.
The result now follows by applying the criterion for vertex accessibility given in
Proposition \ref{vertices}.
\qede \bigskip

The result below follows from \cite{[T]} or \cite{[BFR]}, but we provide a simpler
argument.

\begin{prop}\label{jacobi}
Jacobi matrices are spectrally complete.
\end{prop}
\noindent \textbf{Proof}
The key point is to notice that for a Jacobi matrix $J = Q^T\, \Lambda\, Q$, the matrix
$\Lambda$ has simple spectrum and the first column $q_1$ of $Q$ can be taken to be a 
collection of strictly positive numbers (\cite{[DNT]}). Also, $Q = \GS{M}$, for
$M$ having columns $q_1, \Lambda\, q_1, \ldots, \Lambda^{n-1}\, q_1$.
It is easy now to check that the relevant minors of $Q$ are invertible, by computing
Vandermonde determinants.
\qede \bigskip

Finally, we identify $\bar U$ with the spectral polytope $\PS$. We make use of 
the following simple geometric fact.

\begin{lem}
Let $P$ be a convex polytope in $\Rset^n$ with nonempty interior. Let $\gamma(t)$
be a smooth path passing by a vertex $v$ of $P$ at $t=0$. Then $\gamma'(0) = 0 .$
\end{lem}

\noindent \textbf{Proof} After
an appropriate composition with an affine linear transformation taking $v$ to the origin,
we may suppose that $P$ lies in the positive octant of $\Rset^n$
and $v = 0$. Each coordinate of $\gamma(t)$ then is a smooth function, taking 
nonnegative values and equal to zero at $t=0$. By taking Newton quotients, 
each partial derivative equals 0 at zero.
\qede \bigskip

\begin{prop}\label{politopo}
Let $S$ be symmetric, irreducible. The associated convex polytopes $\bar U$ and $\PS$ are equal.
\end{prop}
\noindent \textbf{Proof}
Every vertex of $\PS$ is the image under $F_{\bar \FS}$ of an accessible vertex 
in $\bar \FS$: in particular, since $F_{\bar \FS}$ is a homeomorphism between $\bar \FS$
and $\bar U$, we must have that $\PS \subset \bar{U}$. Now, suppose $v$ is 
a vertex of $\bar{U}$. There must be a matrix $S_v \in \partial \FS$ for which
$F_{\bar \FS}(S_v) = v$. If $S_v$ is not diagonal, there is a minimal nontrivial 
invariant canonical subspace $V_I$ of dimension at least 2. Thus, from Propositions \ref{quebra}
and \ref{variedade},
the slice ${\mathcal F}_{S_v}$ through $S_v$ is at least one dimensional --- 
said differently, there is some path $t \rho \subset \Rset^n_0$ for which 
$\Phi(t \rho, S_v) \subset {\mathcal F}_{S_v}$ passes by $S_v$
with nonzero derivative at $t= 0$. Also, the image of this path under $F_{\bar \FS}$, 
by Proposition \ref{faces}, 
is a path in $\bar{U}$ passing by $v$ with nonzero derivative for $t=0$. 
But this contradicts the fact that $v$ is a vertex of $\bar{U}$, by the previous lemma. 
Thus, $S_v$ is necessarily a diagonal matrix,
and $v$ then is the image of an accessible vertex. Thus every vertex of $\bar{U}$ is an extremal vertex  and $\bar{U} \subset \PS$.
\qede \bigskip

The proof of Theorem 3 is now complete.

%%%%%%%%%%%%%%%%%%%%%%%%%%%%%%%%%%%%%%%%%%%%%%%%%%%%%%%%%%%%%%%%%%%%%%%%%%%%%%%%%%%%%%%%%%%%%%%%

\bigskip

{\obeylines
R. S. Leite and Carlos Tomei
Depto. de Matem{\'a}tica, PUC-Rio
R. Mq. de S. Vicente 225
Rio de Janeiro, RJ 22453-900, Brazil
rsl@mat.puc-rio.br
tomei@mat.puc-rio.br

\end{document}